\documentclass[12pt,reqno]{amsart}
\usepackage{amsmath,amsthm, amscd, amsfonts, amssymb, graphicx}
\usepackage[bookmarksnumbered, colorlinks, plainpages]{hyperref}

\newtheorem{theorem}{Theorem}[section]
\newtheorem{lemma}[theorem]{Lemma}

\newtheorem{corollary}[theorem]{Corollary}
\theoremstyle{definition}

\theoremstyle{remark}
\newtheorem{remark}[theorem]{Remark}

\numberwithin{equation}{section}

\newcommand{\f}{\frac}
\newcommand{\ov}{\overline}
\newcommand{\eps}{\varepsilon}
\newcommand{\ds}{\displaystyle}
\newcommand{\q}{\quad}

\numberwithin{equation}{section}

\begin{document}

\title[Stability of the first order linear recurrence]{On the stability of the first order linear recurrence in topological vector spaces}
\author[M.S.~Moslehian, D.~Popa]{Mohammad Sal Moslehian$^1$ and Dorian Popa$^2$}

\address{$^1$ Department of Pure Mathematics, Center of Excellence in
Analysis on Algebraic Structures (CEAAS), Ferdowsi University of
Mashhad, P.O. Box 1159, Mashhad 91775, Iran.}
\email{moslehian@ferdowsi.um.ac.ir and moslehian@ams.org}
\urladdr{\url{http://profsite.um.ac.ir/~moslehian/}}

\address{$^2$ Technical University, Department of Mathematics, Str. C. Daicoviciu 15, 400020 Cluj-Napoca, Romania.}
\email{Popa.Dorian@math.utcluj.ro}

\keywords{stability; first order linear recurrence; topological vector spaces; convex hull; balanced hull}

\subjclass[2010]{Primary 39B82; secondary 39A10; 39B72.}

\begin{abstract}
Suppose that $\mathcal{X}$ is a sequentially complete Hausdorff locally convex
space over a scalar field $\mathbb{K}$, $V$ is a bounded subset of $\mathcal{X}$,
$(a_n)_{n\ge 0}$ is a sequence in $\mathbb{K}\setminus\{0\}$ with the
property\ $\ds\liminf_{n\to\infty} |a_n|>1$ and $(b_n)_{n\ge 0}$ is
a sequence in $\mathcal{X}$. We show that for every sequence $(x_n)_{n\ge 0}$
in $\mathcal{X}$ satisfying
\begin{eqnarray*}
x_{n+1}-a_nx_n-b_n\in V\q(n\geq 0)
\end{eqnarray*}
there exists a unique sequence $(y_n)_{n\ge 0}$ satisfying the
recurrence $y_{n+1}=a_ny_n+b_n\,\,(n\geq 0)$ and for every $q$ with
$1<q<\ds\liminf_{n\to\infty} |a_n|$, there exists $n_0\in
\mathbb{N}$ such that
\begin{eqnarray*}
x_n-y_n\in \ds\f{1}{q-1}\ov{conv(V^b)}\q (n\geq n_0).
\end{eqnarray*}
\end{abstract}
\maketitle

%========================================================================================%

\section{Introduction}

The stability problem of functional equations was originally raised
by Ulam \cite{17} in 1940 on a talk at Wisconsin University. The
problem posed by Ulam was the following: ``Under what conditions
does there exist an additive mapping near an approximately additive
mapping?'' The first answer to the question was given by Hyers in
the case of Banach spaces \cite{5}. After Hyers' result many papers
dedicated to this topic extending Ulam's problem to other functional
equations and generalizing Hyers' result in various directions were
published (see e.g. \cite{CZE, FOR, 6, JUN, 7, MOS, 13}). As
mentioned in \cite{1} there are much less results on the stability
for functional equations in single variable than in several
variables. A particular case of equation in single variable is the
linear recurrence (difference equation)
\begin{equation}
\label{1.1}
x_{n+1}=a_nx_n+b_n.
\end{equation}
The results on stability of recurrences play an important role in
the theory of dynamical systems and computer science in connection
to the notions of shadowing and controlled chaos (see e.g. \cite{9,
10, 15}). The first result on Hyers--Ulam stability of the linear
recurrence \eqref{1.1} was given by Popa  \cite{11} in the case of
Banach spaces as follows.
\begin{theorem}\label{theorem} \cite{11}
Let $\mathcal{X}$ be a Banach space over the field $\mathbb{K}$ ($\mathbb{R}$ or $\mathbb{C}$),
$(\eps _n)_{n\ge 0}$ a sequence of positive numbers, $(a_n)_{n\ge 0}$ a sequence
in $\mathbb{K}\setminus\{0\}$ with the property
\begin{equation}
\label{1.2} \ds\limsup_{n\to\infty} \ds\f{\eps _n}{\eps
_{n-1}|a_n|}<1 \q\mbox{or}\q \ds\liminf_{n\to\infty} \ds\f{\eps
_n}{\eps _{n-1}|a_n|}>1
\end{equation}
and $(b_n)_{n\ge 0}$ a sequence in $\mathcal{X}$.

Then there exists a constant $L>0$ such that for every sequence $(x_n)_{n\ge 0}$
in $\mathcal{X}$ satisfying the relation
\begin{equation}
\label{1.3} \|x_{n+1}-a_nx_n-b_n\|\le \eps _n\q(n\geq 0)
\end{equation}
there exists a sequence $(y_n)_{n\ge 0}$ given by the linear recurrence
$$y_{n+1}=a_ny_n+b_n\q(n\geq 0)$$
with the property
\begin{equation}
\label{1.4} \|x_n-y_n\|\le L\, \eps _{n-1}.
\end{equation}
for some $n_0\geq 0$ and all $n \geq n_0$.
\end{theorem}

The above  result was extended later by Brzdek, Popa and Xu to the
linear recurrences of higher order with constant coefficients
\cite{12} and to nonlinear recurrences \cite{2} (see also \cite{3}).

The goal of this paper is to extend Theorem \ref{theorem} to the
stability of the linear recurrence \eqref{1.1} in topological vector
spaces (see also \cite{C-A}).

If $\mathcal{X}$ is a topological vector space over the field $\mathbb{K}$,
$(a_n)_{n\ge 0}$ a sequence in $\mathbb{K}$, $(b_n)_{n\ge 0}$ a sequence in
$\mathcal{X}$ and $V$ a bounded subset of $\mathcal{X}$, then the recurrence \eqref{1.1}
is said to be stable in the spirit of Hyers--Ulam whenever there
exists a bounded subset $W$ of $\mathcal{X}$ such that for any sequence
$(x_n)_{n\ge 0}$ in $\mathcal{X}$ satisfying
\begin{equation}
\label{1.5} x_{n+1}-a_nx_n-b_n\in V\q(n\geq 0)
\end{equation}
there exists a  sequence $(y_n)_{n\ge 0}$
satisfying the linear recurrence
\begin{equation}
\label{1.6} y_{n+1}=a_ny_n+b_n,\q  (n\ge 0)
\end{equation}
and
\begin{equation}
\label{1.7}
x_n-y_n\in W
\end{equation}
for some $n_0\geq 0$ and all $n \geq n_0$. For a subset $A$ of a
topological vector space $\mathcal{X}$ over the field $\mathbb{K}$ we denote the
closedness, the balanced hall and the convex hull of $A$ by
$\ov{A}$, $A^b$ and $convA$, respectively. Recall that for $\lambda
,\mu \in \mathbb{K}$ the following relation holds $ (\lambda +\mu )A\subseteq
\lambda A+\mu A. $ If $\mathbb{K}$ is one of the fields $\mathbb{R}$ or
$\mathbb{C}$, $A$ is a convex set and $\lambda ,\mu $ are positive
numbers, then $ (\lambda +\mu )A=\lambda A+\mu A. $ The reader is
referred to \cite{RUD} for undefined notation and terminology.

%========================================================================================%
\section{Main results}
We give, for the beginning, two auxiliary lemmas which will be used to obtain the main result
of this paper.
Throughout this section $\mathbb{N}_0,\mathbb{R},\mathbb{C}$ stand as usual
for the set of all nonnegative integers, real numbers and complex numbers,
respectively.
By $\mathbb{K}$ we denote one of the fields $\mathbb{R}$ or $\mathbb{C}$.
%========================================================================================%
\begin{lemma}\label{lem1} Let $\mathcal{X}$ be a sequentially complete Hausdorff locally convex space
over $\mathbb{K}$. Suppose that $(a_n)_{n\ge 0}$ is a sequence in
$\mathbb{K}\setminus\{0\}$ such that the series $\ds\sum_{n=0}^\infty a_n$ is
absolutely convergent and let $(v_n)_{n\ge 0}$ be a bounded sequence
in $\mathcal{X}$. Then the series
$$\sum_{n=0}^\infty a_n v_n$$
is convergent in $\mathcal{X}$.
\end{lemma}
\begin{proof} Let $V$ be an arbitrary neighborhood of $0$ and set
$$\sigma _n=\sum_{k=0}^n a_kv_k,\q
s_n=\sum_{k=0}^n |a_k|\q(n\geq 0).$$

We have to prove that there exists $n_V\in \mathbb{N}_0$ such that
$$\sigma _{n+p}-\sigma _n\in V$$
for every $n\ge n_V$ and every $p\in \mathbb{N}_0$.

Let $U$ be a convex and balanced neighborhood of $0$ such that $U\subseteq V$.

Since $A=\{v_n\mid n\in \mathbb{N}\}$ is a bounded set, it follows that
there exists $\alpha >0$ such that $A\subseteq \alpha U$.

From the convergence of $\ds\sum_{n=0}^\infty |a_n|$ follows that there exists
$n_0\in \mathbb{N}$ such that
$$s_{n+p}-s_n<\ds\f{1}{\alpha }\q (n\geq n_0, p\in \mathbb{N}).$$

We have
\begin{align*}
\sigma _{n+p}-\sigma _n
& =\sum_{k=n+1}^{n+p}a_kv_k\\
& =(s_{n+p}-s_n)\sum_{k=n+1}^{n+p}\ds\f{a_k}{s_{n+p}-s_n}\cdot v_k\\
& \in \alpha (s_{n+p}-s_n)\sum_{k=n+1}^{n+p}\ds\f{a_k}{s_{n+p}-s_n}U.
\end{align*}

Since $U$ is a balanced set the following relation holds
$$\ds\f{a_k}{s_{n+p}-s_n}U= \ds\f{|a_k|}{s_{n+p}-s_n}U \q (n+1\le k\le n+p)$$
hence
$$\sigma _{n+p}-\sigma _n\in \alpha (s_{n+p}-s_n)\sum_{k=n+1}^{n+p}\ds\f{|a_k|}{s_{n+p}-s_n}U.$$

The convexity of $U$ leads to
$$\sum_{k=n+1}^{n+p}\ds\f{|a_k|}{s_{n+p}-s_n}U
=\left(\sum_{k=n+1}^{n+p}\ds\f{|a_k|}{s_{n+p}-s_n}\right)U=U$$
therefore
$$\sigma _{n+p}-\sigma _n\in \underbrace{\alpha (s_{n+p}-s_n)}_{<1}U\subseteq U\subseteq V$$
for every $n\ge n_0$ and every $p\in \mathbb{N}$.
\end{proof}
%========================================================================================%
\begin{lemma}\label{lem2} Let $\mathcal{X}$ be a vector space over $\mathbb{K}$, $(a_n)_{n\ge
0}$ be a sequence in $\mathbb{K}$, $(b_n)_{n\ge 0}$ a sequence in $\mathcal{X}$ and
$(x_n)_{n\ge 0}$ a sequence in $\mathcal{X}$ satisfying the recurrence
\eqref{1.1}. Then
\begin{eqnarray*}
x_n=a_0a_1\dots a_{n-1}x_0+\sum_{k=1}^{n-1}a_k\dots
a_{n-1}b_{k-1}+b_{n-1}\q (n\ge 2)\,.
\end{eqnarray*}
\end{lemma}
\begin{proof}
Induction on $n$.
\end{proof}
The first stability result for recurrence \eqref{1.1} is given in
the next theorem.
%========================================================================================%
\begin{theorem} \label{th1} Suppose that $\mathcal{X}$ is a sequentially complete Hausdorff locally
convex space over $\mathbb{K}$, $V$ is a bounded subset of $\mathcal{X}$, $(a_n)_{n\ge
0}$ is a sequence in $\mathbb{K}\setminus\{0\}$ with the property
$\ds\liminf_{n\to\infty} |a_n|>1$ and $(b_n)_{n\ge 0}$ is a sequence
in $\mathcal{X}$. Then for every sequence $(x_n)_{n\ge 0}$ in $\mathcal{X}$ satisfying
\begin{equation}
\label{2.2} x_{n+1}-a_nx_n-b_n\in V\q(n\geq 0)
\end{equation}
there exists a unique sequence $(y_n)_{n\ge 0}$ satisfying the recurrence
\begin{equation}
\label{2.3} y_{n+1}=a_ny_n+b_n\q(n\geq 0)
\end{equation}
and for every $q$ with $1<q<\ds\liminf_{n\to\infty} |a_n|$, there
exists $n_0\in \mathbb{N}$ such that
\begin{equation}
\label{2.4} x_n-y_n\in \ds\f{1}{q-1}\ov{conv(V^b)}\q (n\geq n_0)\,.
\end{equation}
\end{theorem}

\begin{proof} \textbf{Existence.} Let $(x_n)_{n\ge 0}$ be a sequence in $\mathcal{X}$
satisfying (\ref{2.2}) and define the sequence $(c_n)_{n\ge 0}$ by
\begin{equation}
\label{2.5} c_n=x_{n+1}-a_nx_n-b_n\q(n\geq 0).
\end{equation}

Taking account of
\begin{eqnarray*} \ds\limsup_{n\to\infty} \ds\f{\ds\f{1}{|a_0\dots
a_{n+1}|}}{\ds\f{1}{|a_0\dots a_n|}} =\ds\limsup_{n\to\infty}
\ds\f{1}{|a_{n+1}|}<1
\end{eqnarray*}
it follows that the series
\begin{equation}
\label{2.7}
\sum_{n=0}^\infty \ds\f{1}{|a_0\dots a_n|}
\end{equation}
is convergent in view of D'Alembert ratio test.

The boundedness of $(c_n)_{n\ge 0}$ and the convergence of the
series (\ref{2.7}) implies the convergence of the series
$\ds\sum_{n=0}^\infty \ds\f{c_n}{a_0a_1\dots a_n}$ in view of Lemma
\ref{lem1}. Put
\begin{eqnarray*} \sum_{n=0}^\infty \ds\f{c_n}{a_0a_1\dots a_n}=s \q (s\in \mathcal{X})
\end{eqnarray*}
and define the sequence $(y_n)_{n\ge 0}$, by the recurrence
(\ref{2.3}) with
\begin{eqnarray*}
y_0=x_0+s.
\end{eqnarray*}

It follows from Lemma \ref{lem2} that
\begin{align}
\label{2.10}
x_n-y_n
& =-\prod_{k=0}^{n-1}a_ks+\sum_{k=1}^{n-1}a_k\dots a_{n-1}c_{k-1}+c_{n-1}\nonumber\\
& =a_0a_1\dots a_{n-1}\left(-s+\sum_{k=1}^n \ds\f{c_{k-1}}{a_0\dots a_{k-1}}\right)\nonumber\\
& =a_0\dots a_{n-1}\sum_{k=0}^\infty \ds\f{c_{n+k}}{a_0a_1\dots a_{n+k}}\nonumber\\
& =\sum_{k=0}^\infty \ds\f{c_{n+k}}{a_na_{n+1}\dots a_{n+k}}.
\end{align}

Let $q$ be a real number such that
\begin{eqnarray*}
1<q<\ds\liminf_{n\to\infty} |a_n|.
\end{eqnarray*}

Then, there exists $n_0\in \mathbb{N}$ such that
$|a_n|\ge q$ for every $n\ge n_0$.

The following relations hold
\begin{align}
\label{2.12}
\ds\f{c_{n+k}}{a_n\dots a_{n+k}}
& \in \ds\f{1}{a_n\dots a_{n+k}}V\subseteq
\ds\f{1}{a_n\dots a_{n+k}}V^b=\ds\f{1}{|a_n\dots a_{n+k}|}V^b\nonumber\\
& \subseteq \ds\f{1}{q^{k+1}}V^b \subseteq
\ds\f{1}{q^{k+1}}conv(V^b) \q(n,k\in \mathbb{N}_0, n\ge n_0).
\end{align}

From (\ref{2.10}) and (\ref{2.12}) we have
\begin{align}
\label{2.13}
\sum_{k=0}^p \ds\f{c_{n+k}}{a_n\dots a_{n+k}}
& \in \sum_{k=0}^p \ds\f{1}{q^{k+1}}conv(V^b)\nonumber\\
& =\left(\sum_{k=0}^p \ds\f{1}{q^{k+1}}\right)conv(V^b)\q (n\ge
n_0)\,.
\end{align}
By letting $p\to \infty $ in \eqref{2.13} we get
\begin{eqnarray*}
x_n-y_n\in \ds\f{1}{q-1}\ov{conv(V^b)}\q (n\geq n_0).
\end{eqnarray*}

The existence is proved.

\textbf{Uniqueness.} Let $(x_n)_{n\ge 0}$ be a sequence satisfying
(\ref{2.2}) and suppose that there exists a sequence $(y_n)_{n\ge
0}$ satisfying (\ref{2.3}) and (\ref{2.4}) with $y_0\ne x_0+s$. We
have
\begin{align*}
x_n-y_n
& =\prod_{k=0}^{n-1}a_k(x_0-y_0)+\sum_{k=1}^{n-1}a_k\dots a_{n-1}c_{k-1}+c_{n-1}\\
& =\prod_{k=0}^{n-1}a_k\left(x_0-y_0+\sum_{k=1}^n \ds\f{c_{k-1}}{a_0a_1\dots a_{k-1}}\right).
\end{align*}

Since
$$\ds\lim_{n\to \infty }\left(x_0-y_0+\sum_{k=1}^n \ds\f{c_{k-1}}{a_0a_1\dots a_{k-1}}\right)
=x_0-y_0+s\ne 0$$
and
$$\ds\lim_{n\to \infty }\left|\prod_{k=0}^{n-1}a_k\right|=\infty ,$$
in view of the convergence of $\ds\sum_{n=0}^\infty
\ds\f{1}{|a_0\dots a_n|}$, it follows that $(x_n-y_n)_{n\ge 0}$ is
an unbounded sequence, a contradiction to (\ref{2.4}).
\end{proof}
%========================================================================================%

A similar result holds for the linear recurrence with constant coefficients.

\begin{theorem}
\label{t2.4}
Let $\mathcal{X}$ be a sequentially complete Hausdorff locally convex space over $\mathbb{K}$,
$V$ a bounded subset of $\mathcal{X}$, $a\in \mathbb{K}$, $|a|>1$, and $(b_n)_{n\ge 0}$ a sequence in $\mathcal{X}$.
Then for every sequence $(x_n)_{n\ge 0}$ in $\mathcal{X}$ satisfying
\begin{equation}
\label{2.9a} x_{n+1}-ax_n-b_n\in V\q (n\ge 0)
\end{equation}
there exists a unique sequence $(y_n)_{n\ge 0}$ fulfilling the recurrence
\begin{equation}
\label{2.10a} y_{n+1}=ay_n+b_n\q (n\ge 0)
\end{equation}
such that
\begin{eqnarray*}
x_n-y_n\in \ds\f{1}{|a|-1}\cdot \ov{conv(V^b)}\q (n\ge 0).
\end{eqnarray*}
\end{theorem}

\begin{proof}
Denoting $c_n:=x_{n+1}-ax_n-b_n\,\,(n\ge 0)$ it follows that the
series $\ds\sum_{n=0}^\infty \ds\f{c_n}{a^{n+1}}$ is convergent in
view of Lemma \ref{lem1}. Put
$$\sum_{n=0}^\infty \ds\f{c_n}{a^{n+1}}=s \q (s\in \mathcal{X})$$
and define $(y_n)_{n\ge 0}$ by the recurrence (\ref{2.10a}) with $y_0=x_0+s$.
It follows, as in the proof of Theorem \ref{th1}, that
$$x_n-y_n=\sum_{k=0}^\infty \ds\f{c_{n+k}}{a^{k+1}}\q (n\ge 0)$$
and
\begin{eqnarray*}
\ds\f{c_{n+k}}{a^{k+1}}\in \ds\f{1}{|a|^{k+1}}\cdot V^b\subseteq \ds\f{1}{|a|^{k+1}}\cdot conv(V^b).
\end{eqnarray*}
Then
\begin{equation}
\label{2.13a}
\sum_{k=0}^p \ds\f{c_{n+k}}{a^{k+1}}\in \left(\sum_{k=0}^p \ds\f{1}{|a|^{k+1}}\right)conv(V^b).
\end{equation}
Now by letting $p\to \infty $ in (\ref{2.13a}) we get (\ref{2.10a}).

The uniqueness follows analogously to Theorem \ref{th1}.
\end{proof}

\begin{corollary} Suppose that $\mathcal{X}$ is a Banach space over $\mathbb{K}$, $\varepsilon >0$, $|a|>1$ and $(b_n)_{n\ge 0}$ is a sequence in $\mathcal{X}$. Then for every
sequence $(x_n)_{n\ge 0}$ in $\mathcal{X}$ satisfying
\begin{eqnarray*}
\| x_{n+1}-ax_n-b_n\| \leq \eps\q(n\geq 0)
\end{eqnarray*}
there exists a unique sequence $(y_n)_{n\ge 0}$ satisfying the
recurrence
\begin{eqnarray*}
 y_{n+1}=ay_n+b_n\q(n\geq 0)
\end{eqnarray*}
such that
\begin{eqnarray*}
\|x_n-y_n\| \leq \f{\eps}{|a|-1}\q (n\geq 0).
\end{eqnarray*}
\end{corollary}
\begin{proof}
Use Theorem \ref{t2.4} with $a_n=a$ and take $V$ to be the closed
ball of center $0$ with radius $\eps$.
\end{proof}
%========================================================================================%
\begin{remark}
If $\ds\liminf_{n\to\infty}|a_n|\leq 1$, then the conclusion of
Theorem \ref{th1} may be not true.\\ To see this, set
$\mathcal{X}=\mathbb{K}=\mathbb{R}$, take $1\leq r <2$, let $V$ be the interval
$(-r,r)$, $a_n=r, b_n=0\,\,(n\geq 0)$ and consider the sequence
$(x_n)_{n\ge 0}$ given by $x_{n+1}-rx_n=1\,\,(n\geq 0), x_0=0$. Then
$x_n=\sum_{j=1}^{n-1}r^j$. One can observe that for any sequence
$(y_n)_{n\ge 0}$ satisfying the recurrence $y_{n+1}=ry_n$ we have
$$\sup_{n\to\infty}|x_n-y_n|=\infty\,.$$
In fact,
$$\sup_{n\to\infty}|x_n-y_n|=\sup\left(\left\{\left|\sum_{j=0}^{n-1}r^j-r^ny_0\right|: n=1, 2, \cdots\right\}\cup\{|y_0|\}\right)\,.$$
If $y_0\leq 0$, then $\sup_{n\to\infty}|x_n-y_n|=\infty$. Let us
assume $y_0>0$. There exist positive integers $k_0$ and $n_0$ such
that $r^{k_0}\geq y_0$ and $\sum_{j=0}^{n-1}r^j \geq r^n$ for all
$n\geq n_0$. Then
$$\sum_{j=k_0}^{n-1+k_0}r^j-r^ny_0\leq r^{k_0}\left(\sum_{j=0}^{n-1}r^j - r^n\right)=(2-r)r^{n+k_0}-r^{k_0}$$
for all $n\geq \max\{k_0, n_0\}$. Therefore
$\left|\sum_{j=0}^{n-1}r^j-r^ny_0\right|\to \infty$ as $n \to
\infty$. Hence $\sup_{n\to\infty}|x_n-y_n|=\infty$.
\end{remark}
%========================================================================================%

\begin{theorem}\label{th2} Let $\mathcal{X}$ be a Hausdorff topological vector space over the field $\mathbb{K}$,
$V$ a bounded subset of $\mathcal{X}$, $(a_n)_{n\ge 0}$ a sequence in
$\mathbb{K}$ with $\ds\limsup_{n\to\infty} |a_n|<1$ and
$(b_n)_{n\ge 0}$ a sequence in $\mathcal{X}$. Then there exists a positive number $M$ such that for every
sequence $(x_n)_{n\ge 0}$ in $\mathcal{X}$ satisfying
\begin{eqnarray*}
 x_{n+1}-a_nx_n-b_n\in V\q(n\geq 0)
\end{eqnarray*}
there exists a sequence $(y_n)_{n\ge 0}$ satisfying the linear recurrence
\begin{equation}
\label{2.16} y_{n+1}=a_ny_n+b_n\q(n\geq 0)
\end{equation}
such that
\begin{eqnarray*}
x_n-y_n\in M\cdot conv(V^b)\q (n\geq n_0).
\end{eqnarray*}
\end{theorem}
\begin{proof}
Let $(c_n)_{n\ge 0}$ be defined by (\ref{2.5}), as in the proof of Theorem 2.3 and
$(y_n)_{n\ge 0}$ be given by the recurrence (\ref{2.16}) with $y_0=x_0$.
Then
\begin{eqnarray*}
x_n-y_n=\sum_{k=1}^{n-1}a_k\dots a_{n-1}c_{k-1}+c_{n-1} \q (n\ge 2).
\end{eqnarray*}
Choose $q\in \mathbb{R}$ such that
\begin{eqnarray*} \ds\limsup_{n\to\infty} |a_n|<q<1.
\end{eqnarray*}

Then there exists $n_0\in \mathbb{N}$ such that
\begin{eqnarray*}
|a_n|\le q \q (n\geq n_0).
\end{eqnarray*}

For $k\in \mathbb{N}$, $k\ge n_0$ we have
\begin{eqnarray*}
|a_k\dots a_{n-1}|\le q^{n-k}
\end{eqnarray*}
and for $k<n_0$
\begin{align*} |a_k\dots a_{n-1}|
& =|a_k\dots a_{n_0-1}|\cdot |a_{n_0}\dots a_{n-1}|\nonumber\\
& \le q^{n-n_0}|a_k\dots a_{n_0-1}|\nonumber\\
& \le |a_k\dots a_{n_0-1}|.
\end{align*}

Then, the same argument as in the proof of Theorem 2.3 follows
\begin{align*}
x_n-y_n
& =\sum_{k=1}^{n_0-1}a_k\dots a_{n-1}c_{k-1}+\sum_{k=n_0}^{n-1}a_k\dots a_{n-1}c_{k-1}+c_{n-1}\\
& \in \sum_{k=1}^{n_0-1}|a_k\dots a_{n_0-1}|\cdot V^b+\sum_{k=n_0}^{n-1}q^{n-k}\cdot V^b+V^b\\
& \subseteq \sum_{k=1}^{n_0-1}|a_k\dots a_{n_0-1}|conv(V^b)+\ds\f{q}{1-q}conv (V^b)+conv(V^b)\\
& =\left(\sum_{k=1}^{n_0-1}|a_k\dots a_{n_0-1}|+\ds\f{1}{1-q}\right)conv(V^b).
\end{align*}

Choosing
$$M:=\sum_{k=1}^{n_0-1}|a_k\dots a_{n_0-1}|+\ds\f{1}{1-q}\,.$$
\end{proof}
%========================================================================================%
In Theorem \ref{th2} the sequence $(y_n)_{n\ge 0}$ is not necessary uniquely
determined (see \cite[Remark 2.2]{11}).
%========================================================================================%

\begin{theorem}
\label{t2.8} Let $\mathcal{X}$ be a Hausdorff topological vector space over
the field $\mathbb{K}$, $V$ a bounded subset of $\mathcal{X}$, $a\in \mathbb{K}$, $|a|<1$, and
$(b_n)_{n\ge 0}$ a sequence in $\mathcal{X}$. Then for every sequence
$(x_n)_{n\ge 0}$ in $\mathcal{X}$ satisfying \eqref{2.9a} there exists a
sequence $(y_n)_{n\ge 0}$ in $\mathcal{X}$ fulfilling the recurrence
(\ref{2.10a}) such that
$$x_n-y_n\in \ds\f{1}{1-|a|}\cdot conv(V^b) \q (n\ge 0).$$
\end{theorem}

\begin{proof}
Let $c_n:=x_{n+1}-ax_n-b_n\,\,(n\ge 0)$ and $(y_n)_{n\ge 0}$ be
given by the recurrence (\ref{2.10a}) with $y_0=x_0$. Then
$$x_n-y_n=\sum_{k=1}^n a^{n-k}c_{k-1} \q (n\ge 1).$$
It follows that
\begin{align*}
x_n-y_n
& \in \sum_{k=1}^n a^{n-k}V\subseteq \sum_{k=1}^n |a|^{n-k}V^b\\
& \subseteq \sum_{k=1}^n |a|^{n-k}conv(V^b)\\
& =\left(\sum_{k=1}^n |a|^{n-k}\right)conv(V^b)\\
& =\ds\f{1-|a|^n}{1-|a|}conv(V^b)\\
& \subseteq \ds\f{1}{1-|a|}conv(V^b) \q (n\ge 0).
\end{align*}
\end{proof}

\begin{corollary} Suppose that $\mathcal{X}$ is a normed space over $\mathbb{K}$, $\varepsilon >0$, $|a|<1$ and $(b_n)_{n\ge 0}$ is a sequence in $\mathcal{X}$. Then for every
sequence $(x_n)_{n\ge 0}$ in $\mathcal{X}$ satisfying
\begin{eqnarray*}
\| x_{n+1}-ax_n-b_n\| \leq \eps\q(n\geq 0)
\end{eqnarray*}
there exist a positive number $M$ and a sequence $(y_n)_{n\ge 0}$
satisfying the recurrence
\begin{eqnarray*}
 y_{n+1}=ay_n+b_n\q(n\geq 0)
\end{eqnarray*}
such that
\begin{eqnarray*}
\|x_n-y_n\| \leq \ds\f{1}{1-|a|}\eps\q (n\geq 0).
\end{eqnarray*}
\end{corollary}
\begin{proof}
Use Theorem \ref{t2.8} with $a_n=a$ and take $V$ to be the closed
ball of center $0$ with radius $\eps$.
\end{proof}
%========================================================================================%


\begin{thebibliography}{99}
\bibitem{1}
R.P. Agarwal, B. Xu and W. Zhang, \textit{Stability of functional
equations in single variable}, J. Math. Anal. Appl. \textbf{288}
(2003), 852-869.

\bibitem{2}
J. Brzdek, D. Popa and B. Xu, \textit{The Hyers--Ulam stability of
the nonlinear recurrences}, J. Math. Anal. Appl. \textbf{335}
(2007), 443-449.

\bibitem{3}
J. Brzdek, D. Popa and B. Xu, \textit{Note on the nonstability of
the linear recurrence}, Abh. Math. Sem. Univ. Hamburg \textbf{76}
(2006), 183-189.

\bibitem{CZE} S. Czerwik, \textit{Functional Equations and Inequalities in Several Variables},
World Scientific, New Jersey, London, Singapore, Hong Kong, 2002.

\bibitem{C-A} S. Czerwik and M. Adam, \textit{On the stability of the quadratic functional equation in topological
      spaces}, Banach J. Math. Anal. \textbf{1} (2007), no. 2, 245--251.


\bibitem{FOR} G.L.~Forti, \textit{Hyers--Ulam stability of functional equations in
several variables}, Aequationes Math. \textbf{50} (1995), 143--190.

\bibitem{5}
D.H. Hyers, \textit{On the stability of the linear functional
equation}, Proc. Nat. Acad. Sci. USA \textbf{27} (1941), 222-224.

\bibitem{6}
D.H. Hyers, G. Isac and Th.M. Rassias, \textit{Stability of
Functional Equations in Several Variables}, Birkh\"auser, Boston,
1998.

\bibitem{JUN} S.-M. Jung, \textit{Hyers--Ulam--Rassias Stability of Functional Equations in Mathematical Analysis},
Hadronic Press, Palm Harbor, Florida, 2001.

\bibitem{7}
M. Kuczma, \textit{Functional Equations in Single Variable}, Polish
Scientific Publishers, Warszawa, 1968.

\bibitem{MOS} M.S. Moslehian, \textit{Ternary derivations, stability and physical
aspects}, Acta Applicandae Math. \textbf{100} (2008), no. 2,
187--199.

\bibitem{9}
K. Palmer, \textit{Shadowing in Dynamical Systems}, Kluwer Academic
Press, 2000.

\bibitem{10}
S. Pilyugin, \textit{Shadowing in Dynamical Systems}, Lecture Notes
in Mathematics 1706, Springer-Verlag, 1999.

\bibitem{11}
D. Popa, \textit{Hyers--Ulam--Rassias stability of a linear recurrence}, J.
Math. Anal. Appl. \textbf{309} (2005), 591-597.

\bibitem{12}
D. Popa, \textit{Hyers--Ulam stability of the linear recurrence with
constant coefficients}, Adv. in Diff. Eq. 2005-2 (2005), 101-107.

\bibitem{13}
Th.M. Rassias(ed.), \textit{Functional Equations, Inequalities and
Applications}, Kluwer Academic Publishers, Dordrecht, Boston and
London, 2003.

\bibitem{RUD}
W. Rudin, \textit{Functional analysis}, Second edition.
International Series in Pure and Applied Mathematics. McGraw-Hill,
Inc., New York, 1991.

\bibitem{15}
J. Tabor, \textit{General stability of functional equations of
linear type}, J. Math. Anal. Appl. \textbf{1} (2007), 192-200.

\bibitem{17}
S.M. Ulam, \textit{Problems in Modern Mathematics}, Wiley, New York,
1964.


\end{thebibliography}
\end{document}